\documentclass[12pt, reqno]{amsart}

\usepackage{amsmath, amscd, amsfonts, amssymb, graphicx, color}

\usepackage[utf8]{inputenc}

\IfFileExists{mathrsfs.sty}{\usepackage{mathrsfs}}%
{\typeout{Package mathrsfs not found: mathscr ---> mathcal}
\newcommand{\mathscr}{\mathcal}}

\usepackage{xspace}
\usepackage{xfrac}
\usepackage{faktor}

\numberwithin{equation}{section}

\newcommand{\BBA}{{\mathbf A}}
\newcommand{\BBC}{{\mathbb C}}
\newcommand{\BBD}{{\mathbf D}}
\newcommand{\BBN}{{\mathbb N}}

\newenvironment{mmatrix}{%
\setlength{\arraycolsep}{3pt}%
\begin{pmatrix}}{\end{pmatrix}}

\usepackage{relsize}

\DeclareMathOperator{\ran}{ran}

\newcommand\quotient[2]{
\text{\raise0.5ex\hbox{$#1$}{\Big/}
\lower0.5ex\hbox{$#2$}}
}

\newcommand{\harp}{\hspace{-0.7ex}\upharpoonright}
\newcommand{\vertic}{\hspace{0.25ex}|\hspace{0.25ex}}

\newcommand{\QED}{$\square$}

\newtheorem{theorem}{Theorem}[section]
\newtheorem{lemma}[theorem]{Lemma}
\newtheorem{corollary}[theorem]{Corollary}
\newtheorem{proposition}[theorem]{Proposition}

\theoremstyle{definition}
\newtheorem{definition}[theorem]{Definition}

\newenvironment{pf}[1][Proof]
{\begin{proof}[#1]}
{\end{proof}}

\AtBeginDocument{%
   \def\MR#1{}
}

\begin{document}

\setcounter{page}{1}

\title[A generalized Schur complement]{A generalized Schur complement
    for non-negative operators on linear spaces}

\author[J.~Friedrich, M.~G\"unther, \MakeLowercase{and} L.~Klotz]%
{J.~Friedrich,$^1$ M.~G\"unther,$^2$$^{\ast}$ \MakeLowercase{and}
  L.~Klotz$^2$}

\address{$^{1}$Stauffenbergstr. 10, D-04509 Delitzsch, Germany.}
\email{\textcolor[rgb]{0.00,0.00,0.84}{jf.dz@alice.de}}

\address{$^{2}$Mathematisches Institut, Universit\"at Leipzig, 
PF 10\,09\,20, D-04009 Leipzig, Germany.}
\email{\textcolor[rgb]{0.00,0.00,0.84}{guenther@math.uni-leipzig.de;
    klotz@math.uni-leipzig.de}}



\date{\today
\newline\indent $^{*}$Corresponding author}

\begin{abstract}
Extending the corresponding notion for matrices or bounded
linear operators on a Hilbert space we define a generalized Schur
complement for a non-negative linear operator mapping a linear space
into its dual and derive some of its properties.
\end{abstract}

\maketitle



\section{Introduction}
\label{intro}
In case of $2\times 2$~block matrices the notion of Schur
complement and its generalizations have a long history. We refer to
\cite{Zhang}, where it is given a comprehensive exposition of the
history, theory and versatility of Schur complements. In the case of
non-negative matrices the notion of generalized Schur complement can
be extended to matrices whose entries are bounded operators in Hilbert
spaces, cf.~\cite{Shmulyan}. Moreover, generalized Schur complements are
closely related to shorted operators which were introduced by
M.G.~Kre\u{\i}n \cite{Krein} and have found interesting applications in
electrical network theory, cf.~\cite{Anderson}.
Since there is a natural way to define non-negativity for a linear
operator, which maps a linear space into its dual, one can ask for a
generalized Schur complement of such an operator. A first attempt was
made in \cite{FKKlotz}, where non-negative bounded linear operators
from a Banach space into its topological dual were discussed. The
present paper deals with a generalized Schur complement and a shorted
operator of a non-negative $2\times 2$ block matrix, whose entries are
linear operators on linear spaces. Thus, topological questions,
particularly continuity problems, play a minor role.
To define a generalized Schur complement of a non-negative operator on
a linear space, one needs the notion of a square
root. Section~\ref{sec2} deals with this important and useful concept,
which was studied by many authors.
Section~\ref{sec3} contains definitions and basic properties of the Schur
complement and the shorted operator for slightly more general
than non-negative operators. In Section~\ref{sec4} further results on
generalized Schur complements  are derived.
Among other things we extend the Crabtree-Haynsworth quotient formula
\cite{Crabtree}. One of the most useful results concerning $2\times
2$ block matrices is Albert's non-negativity criterion \cite{Albert}.
A generalization to non-negative operators on linear spaces and some of
its consequences are given in Section~\ref{sec5}. The special class
of extremal operators, which was introduced by M.G.~Kre\u{\i}n \cite{Krein}
is the subject of Section~\ref{sec6}.

For bounded linear operators on Hilbert spaces many results concerning
the generalized Schur complement were obtained by Yu. L.~Shmulyan. A
large part of them was proved independently or rediscovered later on
by several mathematicians from western countries. The present paper is
strongly influenced by Shmulyan's work and was written to illustrate
his contribution to the theory of generalized Schur complement. In
this way most of our assertions of Sections~\ref{sec3}-\ref{sec5}
are generalizations of results contained in \cite{Shmulyan} to
non-negative operators on linear spaces.

\section{Basic definitions and notations}
\label{sec1}
Any linear space of the present paper is a space over $\BBC$, the field of
complex numbers, and its zero element is denoted by $0$. For a linear
space $X$, let $X'$ denote its dual space of all antilinear
functionals on $X$ and $\langle x',x\rangle_{X}:=\langle x',x\rangle$
the value of $x'\in X'$ at $x\in X$. If $X^{\sim}$ is a subspace of
$X'$, an arbitrary $x\in X$ defines an element $jx$ of $(X^{\sim})'$
according to
\[
\langle jx,x^{\sim}\rangle_{X'}
:=\overline{\langle x^{\sim},x\rangle}_{X}, \quad x^{\sim}\in
X^{\sim},
\]
where $\bar{\alpha}$ stands for the complex conjugate of $\alpha\in\BBC$.
\par\medskip\noindent
\textbf{Convention (CN):}
If for all $x\in X\setminus\{0\}$ there exists $x^{\sim}\in X^{\sim}$
such that $\langle x^{\sim},x\rangle\neq0$, we shall identify $X$ and
its isomorphic image under the map $j$ and write
\[
  \langle jx,x^{\sim}\rangle_{X'} =: \langle x, x^{\sim}\rangle_{X'},
  \quad x\in X, x^{\sim}\in X^{\sim}.
\]
\par\medskip\noindent
The linear space of all linear operators from $X$ into a linear space
$Y$ is denoted by $ \mathscr L (X,Y)$, and $I$ is the identity
operator in case $X=Y$. If $A\in \mathscr L (X,Y)$ and $X_{1}$ is a
subspace of $X$, the symbols $\ker A$, $\ran A$ and
$A\harp_{X_{1}}$ stand for the
null space, range, and restriction of $A$ to $X_{1}$, resp.
Set $AX_{1}:= \ran A\harp_{X_{1}}$.
The dual operator $A'\in \mathscr L (Y',X')$
is defined by the relation $\langle y',Ax\rangle_{Y}=\langle
A'y',x\rangle_{X}$, $x\in X$, $y'\in Y'$.
\par\medskip\noindent
\textbf{Examples.}
1. If $Z$ is a linear space and $A\in\mathscr L(X,Y)$,
$B\in\mathscr L(Y,Z)$, then $(BA)'=A'B'$. \\
2. If $A\in\mathscr L(X,Y)$, then $A''\in\mathscr L(X'',Y'')$
and $A=A''\harp_{X}$ according to (CN). \\
3. If $A\in\mathscr L(X,X')$, then $A'\in\mathscr L(X'',X')$.
Taking into account (CN), we get
$\langle x_{2},Ax_{1}\rangle_{X'} =
\overline{\langle Ax_{1},x_{2}\rangle}_{X}$ and
$\langle x_{2},Ax_{1}\rangle_{X'} =\langle A'x_{2},x_{1}\rangle_{X}$,
hence,
\begin{equation}
\label{sec1:1}
\langle Ax_{1},x_{2}\rangle_{X} =
\overline{\langle A'x_{2},x_{1}\rangle}_{X},\quad x_{1},x_{2}\in X.
\end{equation}
\par\medskip\noindent
An operator $A\in\mathscr L(X,X')$ is called Hermitian, if
$\langle Ax_{1},x_{2}\rangle =
\overline{\langle Ax_{2},x_{1}\rangle}$ and non-negative if
$\langle Ax_{1},x_{1}\rangle \geq0$, $x_{1},x_{2}\in X$. The sets of
all
Hermitian and all non-negative operators are denoted by
$\mathscr L^{h}(X,X')$ and $\mathscr L^{\geq}(X,X')$, resp..
The polarization identity implies that $A$ is Hermitian if and only if
$\langle Ax,x\rangle$ is real for all $x\in X$. Thus
$\mathscr L^{\geq}(X,X')\subseteq \mathscr L^{h}(X,X')$ and the space
$\mathscr L^{h}(X,X')$ can be provided with Loewner's semi-ordering,
i.e. for $A,D\in\mathscr L^{h}(X,X')$ we shall write $A\leq D$ if and only
if $\langle Ax,x\rangle\leq\langle Dx,x\rangle$, $x\in X$. 
Recall Cauchy's inequality
\begin{equation}
\label{sec1:2}
|\langle Ax_{1},x_{2}\rangle|^{2} \leq
\langle Ax_{1},x_{1}\rangle \langle Ax_{2},x_{2}\rangle,\quad
x_{1},x_{2} \in X,
\end{equation}
if $A\in\mathscr L^{\geq}(X,X')$.

\section{Square roots}
\label{sec2}
Let $H$ be a complex Hilbert space with norm
$\|\cdot\| := \|\cdot\|_{H}$
and inner product $(\cdot\vertic\cdot):=(\cdot\vertic\cdot)_{H}$,
which is assumed to be antilinear with respect to the second
component. Let $R\in \mathscr L(X,H)$. Identifying $H$ and the space of
continuous antilinear functionals on $H$ in the common way, one
has $H\subseteq H'$ and
\begin{equation}
\label{sec2:1}
(h\vertic Rx)=\langle R'h,x\rangle,\quad x\in X, h\in H.
\end{equation}
Set $R^{\ast}:= R'\harp_{H}$. From \eqref{sec2:1} it can be concluded
that $\ker R^{\ast}$ is equal to the orthogonal complement of $(\ran
R)^{c}$, where $M^{c}$ denotes the closure of a subset $M$ of a
topological space. It follows that $R^{\ast}$ is one-to-one if and
only if $\ran R$ is dense in $H$ and that
\begin{equation}
\label{sec2:2}
\ran R^{\ast} = R^{\ast}(\ran R)^{c}
\end{equation}
Therefore, we can define a generalized inverse $R^{\ast[-1]}$
of $R^{\ast}$ by
\[
R^{\ast[-1]} x' := (R^{\ast}\harp_{(\ran R)^{c}})^{-1} x', \quad
x'\in\ran R^{\ast}.
\]
\begin{lemma}
\label{lemma2.1}
Let $R\in\mathscr L(X,H)$. An element $x'\in X'$ belongs to $\ran R^{\ast}$
if and only if the following conditions are satisfied:
\begin{enumerate}
\item[(i)]
If $x\in\ker R$, then $\langle x',x\rangle = 0$.
\item[(ii)]
$\displaystyle\sup_{x\in X}
\frac{|\langle x',x\rangle|^{2}}{\|Rx\|_{H}^{2}} <\infty$ 
(with convention $\frac{0}{0}:=0$ at the left-hand side).
\end{enumerate}
\end{lemma}
\begin{pf}
If $x'\in R^{\ast}h$ for some $h\in H$, then
\[
  |\langle x',x\rangle|=|\langle R^{\ast}h,x\rangle|=|(h\vertic Rx)|
  \leq\|h\|\|Rx\|,
\]
which yields (i) and (ii). Conversely, assume that (i) and (ii)
are satisfied for some $x'\in X'$.
Set $\varphi(Rx):=\langle x',x\rangle$, $x\in X$.
Because of (i) $\varphi$ is correctly defined and (ii) implies that
$\varphi$ is continuous, so that $\varphi$ is a continuous antilinear
functional on $\ran R$. Thus, there exists $h\in H$ such that
$\langle x',x\rangle= (h\vertic Rx)= \langle R^{\ast}h,x\rangle$
for all $x\in X$, which yields $x'=R^{\ast}h\in\ran R^{\ast}$.
\hfill\QED
\end{pf}
\begin{definition}
Let $A\in\mathscr L(X,X')$. A pair $(R,H)$ of a Hilbert space $H$ and
an operator $R\in\mathscr L(X,H)$ is called a square root of $A$ if
$A=R^{\ast}R$,
and a minimal square root if, addionally, $\ran R$ is dense in $H$.
\end{definition}
Note that there exists a square root of $A$ if and only if there
exists a minimal one. The following result is basic to our
considerations and generalizes the fact concerning the existence of a
square root of a non-negative selfadjoint operator in a Hilbert
space. Its well known short proof is recapitulated for convenience of
the reader.
\begin{theorem}
\label{thm2.3}  
An operator $A\in\mathscr L(X,X')$ possesses a square root if and only
if it is non-negative.  
\end{theorem}
\begin{pf}
Let $A\in\mathscr L^{\geq}(X,X')$. Cauchy's inequality \eqref{sec1:2}
implies that
\[
N:=\{x\in X\;:\;\langle Ax,x\rangle =0\}
\]
is a subspace of $X$. Define an inner product on the quotient
space $X/N$ by
\[
  (x_{1}+N\vertic x_{2}+N):=\langle Ax_{1},x_{2}\rangle,
  \quad x_{1},x_{2}\in X,
\]
and denote the completion of the corresponding
inner product space by $H$. Set $Rx:=x+N$, $x\in X$. It follows
$R\in\mathscr L(X,H)$, $(\ran R)^{c}=H$, and
\[
\langle Ax_{1},x_{2}\rangle =(Rx_{1}\vertic Rx_{2})=\langle
R^{\ast}Rx_{1},x_{2}\rangle, \quad x_{1},x_{2}\in X.
\]
Therefore, $(R,H)$ is a minimal square root of $A$. The "only if"-part
of the assertion is obvious.
\hfill\QED
\end{pf}
The notion of a square root of a non-negative operator acting between
spaces more general than Hilbert spaces was discussed and applied by
many authors. Most of them deal with a topological space $X$ and in
this case continuity problems arise additionally. Some properties of
square roots for operators of special type were obtained by Va\u{\i}nberg
and Engel'son, cf. \cite{Vainberg}. For a Banach space $X$, the
construction of the proof of the preceding theorem was published as
an appendix to \cite{Vakhaniya} and attributed to Chobanyan, see also
\cite{Schwartz} and \cite{Vakhaniya1}. Another but related construction
was proposed by Sebesty\'en \cite{Sebestyen},
cf. \cite{Tarcsay}. G\'orniak \cite{Gorniak84} and G\'orniak and Weron
\cite{Gorniak80} dealt with the existence of a continuous square root
if $X$ is a topological linear space. G\'orniak, Makagon and Weron
\cite{Gorniak83} investigated square roots of non-negative
operator-valued measures. Pusz and Woronowicz \cite{Pusz} extended the
construction of the proof of Theorem~\ref{thm2.3} to pairs of non-negative
sequilinear forms, cf. \cite{Titkos} for further generalizations.
\begin{lemma}
\label{lemma2.4}
If $A\in\mathscr L(X,X')$ and $(R,H)$ is a square root of $A$, then
\[
\ker R = \ker A =\{x\in X\;{:}\;\langle Ax,x\rangle = 0\}.
\]
\end{lemma}
\begin{pf}
The result follows from a chain of conclusions:
\[
\langle Ax,x\rangle =0 \Rightarrow \langle R^{\ast}Rx,x\rangle=0
\Rightarrow(Rx\vertic Rx)=0\Rightarrow Rx=0,
\]
and conversely
\begin{equation}
Rx=0 \Rightarrow R^{\ast}Rx=0\Rightarrow
Ax=0\Rightarrow\langle Ax,x\rangle=0.
\tag*{\QED}
\end{equation}
\end{pf}
The preceding results can be used to derive a version of a part of
Douglas' theorem \cite{Douglas}, cf. \cite{Shmulyan1}.
\begin{proposition}
\label{prop2.5}
Let $A,D\in\mathscr L^{\geq}(X,X')$ and $(R_{A},H_{A})$ and
$(R_{D},H_{D})$ be square roots of $A$ and $D$, resp.. The following
assertions are equivalent:
\begin{enumerate}
\item[(i)]
    $A\leq\alpha^{2}D$ for some $\alpha\in[0,\infty)$,
\item[(ii)] there exists a bounded operator $W\in \mathscr L(H_{A},H_{D})$ with operator norm
    $\|W\|\leq\alpha$ and such that $R^{\ast}_{A}=R^{\ast}_{D}W$.
\end{enumerate}
If (i) or (ii) are satisfied, there exists a unique $W$ so
that $W\subseteq(\ran R_{D})^{c}$. Moreover, $\ker W = \ker
R^{\ast}_{A}$ for this operator $W$.
\end{proposition}
\begin{pf}
Since $R^{\ast}_{A}=R^{\ast}_{D}W$ yields
$R_{A}={R^{\ast}_{A}}'\harp_{X}= W'{R^{\ast}_{D}}'\harp_{X} =
W^{\ast}R_{D}$ by (CN), from (ii) it follows
\[
\langle Ax,x\rangle = \|R_{A}x\|^{2} = \|W^{\ast}R_{D}x\|^{2}
\leq \alpha^{2} \|R_{D}x\|^{2}=\alpha^{2}\langle Dx,x\rangle,\quad x\in X,
\]
hence, (i).
Let $W_{j}\in\mathscr L(H_{A},H_{D}$ be such that
$R^{\ast}_{A}=R^{\ast}_{D}W_{j}$
and $\ran W_{j}\subseteq(\ran R_{D})^{c}$, $j=1,2$. Then
$\ran (W_{1}-W_{2})\subseteq \ker R^{\ast}_{D}$ and $\ran
(W_{1}-W_{2})\subseteq(\ran R_{D})^{c}$, which shows that $W_{1}=W_{2}$.
Now assume that (i) is true. One has $\ker R_{D}\subseteq\ker R_{A}$
by Lemma~\ref{lemma2.4}, hence, $\ran R^{\ast}_{A}\subseteq
\ran R^{\ast}_{D}$ by Lemma~\ref{lemma2.1}. The operator
$W:= R^{\ast[-1]}_{D}R^{\ast}_{A}\in\mathscr L(H_{A},H_{D})$ satisfies
$R^{\ast}_{D}W=R^{\ast}_{A}$, $\ker W =\ker R^{\ast}_{A}$ and $\ran
W\subseteq(\ran R_{D})^{c}$.
The inclusion $\ran R^{\ast[-1]}_{D}\subseteq(\ran R_{D})^{c}$ implies
that $W^{\ast}h= R^{\ast\prime}_{A}(R^{\ast[-1]}_{D})^{\ast}h=0$ if $h$ is
orthogonal to $\ran R_{D}$. Therefore, from
\[
\|W^{\ast}R_{D}x\|^{2} = \|R_{A}x\|^{2} = \langle Ax,x\rangle
\leq\alpha^{2} \langle Dx,x\rangle = \alpha^{2} \|R_{D}x\|^{2},
\quad x\in X,
\]
one can conclude that $\|W\|=\|W^{\ast}\|\leq\alpha$.
\hfill\QED
\end{pf}

As a by-product of Proposition~\ref{prop2.5} we obtain the following corollary.

\begin{corollary}
\label{cor2.6}
Let $A$, $D$, $(R_{A},H_{A})$, and $(R_{D},H_{D})$ be as in
Proposition~\ref{prop2.5}.
\begin{enumerate}
\item[(i)]
  If $A\leq\alpha^{2}D$ for some $\alpha\in[0,\infty)$, then $\ran
  R^{\ast}_{A}\subseteq \ran R^{\ast}_{D}$.
\item[(ii)]
  If $\beta^{2}D\leq A\leq \alpha^{2}D$ for some
  $\alpha,\beta\in(0,\infty)$, then $\ran
  R^{\ast}_{A} =\ran R^{\ast}_{D}$.
\item[(iii)]
  If $(S_{A},G_{A})$ is a square root of $A$, then $\ran
  R^{\ast}_{A} =\ran S^{\ast}_{A}$.
\end{enumerate}
\end{corollary}

\begin{corollary}
\label{cor2.7}
Let $H_{j}$ be Hilbert spaces and $R_{j}\in\mathscr L(X,H_{j})$,
$j=1,2$. If $(R,H)$ is a square root of the non-negative operator
$A:=R^{\ast}_{1}R_{1}+R^{\ast}_{2}R_{2}$, then
$\ran R^{\ast}=\ran R^{\ast}_{1}+\ran R^{\ast}_{2}$.
\end{corollary}
\begin{pf}
Let $G$ be the orthogonal sum of $H_{1}$ and $H_{2}$ and
$S\in\mathscr L(X,G)$ be defined by
$S=\begin{mmatrix}R_{1}\\R_{2}\end{mmatrix}$. 
Since $S^{\ast}=(R^{\ast}_{1},R^{\ast}_{2})$ and $S^{\ast}S = A$, we
get that $\ran S^{\ast}=\ran R^{\ast}_{1}+\ran R^{\ast}_{2}$ and that
$(S,G)$ is a square root of $A$. Now apply Corollary~\ref{cor2.6}
(iii).
\hfill\QED
\end{pf}

\begin{lemma}
\label{lemma2.8}
Let $(R,H)$ be a square root and $(S,G)$ a minimal square root
of $A\in\mathscr L^{\geq}(X,X')$. There exists an isometry
$U\in \mathscr L(G,H)$ such that $US=R$.
\end{lemma}

\begin{pf}
By Lemma~\ref{lemma2.4} there exists an operator $\tilde{U}$
satisfying $\tilde{U}Sx=Rx$, $x\in X$. From
$\|\tilde{U}Sx\|^{2}=\|Rx\|^{2} =\langle Ax,x\rangle=\|Sx\|^{2}$
it follows that $\tilde{U}$ is isometric and can be extended to an
isometry $U\in\mathscr L(G,H)$.
\hfill\QED
\end{pf}

\section{Generalized Schur complements and shorted operators of
  operators of positive type}
\label{sec3}
Let $X$ und $Y$ be linear spaces.

\begin{definition}
\label{def3.1}
A pair $(A,B)$ of an operator $A\in\mathscr L^{\geq}(X,X')$ and 
$B\in\mathscr L(Y,X')$ is called a positive pair if 
$\ran B\subseteq \ran R^{\ast}$ for some square root (and, hence, for
all square roots) $(R,H)$ of $A$.
\end{definition}
The following criterion is an immediate consequence of Lemma~\ref{lemma2.1}. 
\begin{lemma}
\label{lemma3.2}
Let $A\in\mathscr L^{\geq}(X,X')$ and $B\in\mathscr L(Y,X')$. The pair $(A,B)$ is 
a positive pair if and only if for all $y\in Y$ the following conditions 
are satisfied:
\begin{enumerate}
\item[(i)]
If $x\in\ker A$, then $\langle By,x\rangle=0$.
\item[(ii)]
$\displaystyle\sup_{x\in X}
\frac{|\langle By,x\rangle|^{2}}{\langle Ax,x\rangle}<\infty$
(with convention $\frac{0}{0}:=0$).
\end{enumerate}
\end{lemma}
According to (CN), the space $X$ can be considered as a subspace of the domain 
of $B'$. To abbreviate the notation we set 
\[
B^{\sim}:= B'\harp_{X}.
\]
Note that $\langle By,x\rangle = \overline{\langle B^{\sim}x,y\rangle}$,
$x\in X$, $y\in Y$. Thus condition (i) of Lemma~\ref{lemma3.2} is 
equivalent to the inclusion $\ker A\subseteq\ker B^{\sim}$.

If $(A,B)$ is a positive pair and $(R,H)$ is a square root of $A$, the operators
\begin{equation}
\label{sec3:1}
T:= R^{\ast[-1]}B
\end{equation}
and
\begin{equation}
\label{sec3:2}
\omega(A,B):=T^{\ast}T = (R^{\ast[-1])}B)^{\ast} R^{\ast[-1]}B
\end{equation}
can be defined. Note that $B=R^{\ast}T$. The following lemma is obvious.
\begin{lemma}
\label{lemma3.3}
If $(A,B)$ is a positive pair and $(R,H)$ is a square root of $A$,
then for all $y\in Y$
\[
\inf\{\|Ty-Rx\|\::\:x\in X\} = 0.
\]
Equivalently, $\ran T\subseteq (\ran R)^{c}$.
\hfill\QED
\end{lemma}
Recall that the dual space of $X\times Y$ can be written as a 
Cartesian product $(X\times Y)'=X'\times Y'$, where 
$\langle\cdot,\cdot\rangle_{X\times Y}=\langle\cdot,\cdot\rangle_{X}+
\langle\cdot,\cdot\rangle_{Y}$. Also, it should not cause confusion if
we identify the subspace $X\times\{0\}$ of  $X\times Y$ with $X$. An 
operator $\BBA$ of $\mathscr L\bigl(X\times Y,(X\times Y)'\bigl)$ can be
represented as a $2\times 2$ matrix 
$\BBA=\begin{mmatrix}A & B \\ C & D\end{mmatrix}$, where 
$A\in \mathscr L(X,X')$, $B,C\in \mathscr L (Y,X')$, 
$D\in \mathscr L(Y,Y')$. It is not hard to see that $\BBA$ is
Hermitian if and only if $A$ and $D$ are Hermitian and $C=B^{\sim}$. 
To abbreviate the notation we set 
$\mathscr L^{h}\bigl(X\times Y,(X\times Y)'\bigr) =: \mathscr L^{h}$ and
$\mathscr L^{\geq}\bigl(X\times Y,(X\times Y)'\bigr) =: \mathscr L^{\geq}$.
\begin{definition}
\label{def3.4}
An operator $\begin{mmatrix}A & B \\ B^{\sim} & D\end{mmatrix}
\in\mathscr L^{h}$ is called an operator of positive type if
$(A,B)$ is a positive pair. The set of operators of positive
type is denoted by 
$\mathscr L^{+}\bigr(X\times Y,(X\times Y)'\bigr) =: \mathscr L^{+}$.
\end{definition}
\begin{definition}
\label{def3.5}
Let $\BBA =\begin{mmatrix}A & B \\ B^{\sim} & D\end{mmatrix}
\in\mathscr L^{+}$. The operator $\sigma(\BBA):=D-\omega(A,B)$ is
called a generalized Schur complement of $\BBA$ and the operator
\[
\mathscr S(\BBA):=
\begin{pmatrix}0 & 0 \\ 0 & \sigma(\BBA)\end{pmatrix}
\]
is called a shorted operator.
\end{definition}
The following result is a generalization of 
\cite[Corollary~1 to Theorem~3]{Anderson}.
\begin{proposition}
\label{prop3.6}
If $\BBA = \begin{mmatrix}A & B \\ B^{\sim} & D\end{mmatrix}
\in\mathscr L^{+}$, then
$\ran \BBA\cap Y' \subseteq \ran\mathscr S(\BBA)$.
\end{proposition}
\begin{pf}
Let $y'\in Y'$ be such that
$\BBA \begin{mmatrix}x \\ y\end{mmatrix} = 
\begin{mmatrix}0 \\ y'\end{mmatrix}$
for some 
$\begin{mmatrix}x \\ y\end{mmatrix} \in
X\times Y$. Let $(R,H)$ be a minimal square root of $A$. 
Since
\[
\BBA=
\begin{pmatrix} R^{\ast}R & R^{\ast}T\\T^{\ast}R&T^{\ast}T
\end{pmatrix} + \mathscr S(\BBA),
\]
one has 
$R^{\ast}Rx+R^{\ast}Ty=0$, hence, $Rx+Ty=0$ and 
$T^{\ast}Rx+T^{\ast}Ty=0$, which yields 
\begin{equation}
\tag*{\QED}
\begin{pmatrix}0 \\ y'\end{pmatrix}
=\mathscr S(\BBA)
\begin{pmatrix}x \\ y\end{pmatrix}
\in \ran\mathscr S(\BBA).
\end{equation}
\end{pf}
The next result is a simple but useful consequence of Lemma~\ref{lemma3.3}.
\begin{proposition}
\label{prop3.7}
Let $\BBA = \begin{mmatrix} A & B \\ B^{\sim} & D \end{mmatrix}
\in\mathscr L^{+}$. For 
$\begin{mmatrix} x \\ y \end{mmatrix}\in X\times Y$,
\begin{equation}
\label{sec3:3}
\left\langle\mathscr S(\BBA)
\begin{pmatrix} x \\ y \end{pmatrix},
\begin{pmatrix} x \\ y \end{pmatrix}
\right\rangle =
\inf_{z\in X} \left\langle
\BBA\begin{pmatrix} x-z \\ y \end{pmatrix}
\begin{pmatrix} x-z \\ y \end{pmatrix}
\right\rangle.
\end{equation}
Particularly, $\sigma(\BBA)$ and $\mathscr S(\BBA)$ do not depend on the choice of the square root of $A$.

\end{proposition}
\begin{pf}
Since \eqref{sec3:3} is independent of $x\in X$, it is enough 
to prove it for $x=0$. From Lemma~\ref{lemma3.3} it follows 
\begin{align}
\notag
& \inf_{z\in X}
\left\langle\BBA
\begin{pmatrix} -z \\ y \end{pmatrix},
\begin{pmatrix} -z \\ y \end{pmatrix}
  \right\rangle \\
  \notag
& \phantom{\inf_{z\in X}}
= \inf_{z\in X}
\left\langle
\begin{pmatrix}
R^{\ast}R & R^{\ast}T \\ T^{\ast}R & T^{\ast}T
\end{pmatrix}
\begin{pmatrix} -z \\ y \end{pmatrix},
\begin{pmatrix} -z \\ y \end{pmatrix}
\right\rangle
+ 
\left\langle\mathscr S(\BBA)
\begin{pmatrix} 0 \\ y \end{pmatrix},
\begin{pmatrix} 0 \\ y \end{pmatrix}
\right\rangle \\
\notag
& \phantom{\inf_{z\in X}}  
= \inf_{z\in X}\|Ty-Rz\|^{2} + 
\left\langle\mathscr S(\BBA)
\begin{pmatrix} 0 \\ y \end{pmatrix},
\begin{pmatrix} 0 \\ y \end{pmatrix}
\right\rangle \\
\tag*{\QED}
& \phantom{\inf_{z\in X}}
=
\left\langle\mathscr S(\BBA)
\begin{pmatrix} 0 \\ y \end{pmatrix},
\begin{pmatrix} 0 \\ y \end{pmatrix}
\right\rangle.
\end{align}
\end{pf}
\begin{corollary}
\label{cor3.8}
(i) If $\BBA\in\mathscr L^{+}$, then $\mathscr S(\BBA)\leq\BBA$
and $\ker\BBA\subseteq\ker\mathscr S(\BBA)$. \\
(ii) If $\BBA, \BBA_{1}\in\mathscr L^{+}$ and $\BBA\leq\BBA_{1}$, then 
$\mathscr S(\BBA)\leq\mathscr S(\BBA_{1})$.
\end{corollary}
\begin{pf}
The first assertion of (i) as well as (ii) are immediately clear from
Proposition~\ref{prop3.7}. To prove the second assertion of (i), let 
$\begin{mmatrix} x \\ y \end{mmatrix}\in\ker\BBA$. If $z\in X$, one has 
\[
\left\langle \BBA
\begin{pmatrix} x-z \\ y \end{pmatrix},
\begin{pmatrix} x-z \\ y \end{pmatrix}
\right\rangle =\langle Az,z\rangle\geq0,
\]
which implies that the infimum at the right hand side of \eqref{sec3:3}
is equal to $0$. Since $\mathscr S(\BBA)\leq \BBA$ and
\[
\left\langle\bigl(\BBA-\mathscr S(\BBA)\bigr)
\begin{pmatrix} x \\ y \end{pmatrix},
\begin{pmatrix} x \\ y \end{pmatrix}
\right\rangle=0,
\]
it follows $\begin{mmatrix} x \\ y \end{mmatrix}
\in\ker\bigl(\BBA-\mathscr S(\BBA)\bigr)$ by Lemma~\ref{lemma2.4},
hence, 
$\begin{mmatrix} x \\ y \end{mmatrix}\in\ker\mathscr S(\BBA)$.
\hfill\QED
\end{pf}

\begin{corollary}
\label{cor3.9}
If $\BBA\in\mathscr L^{\geq}$, then 
$\mathscr S(\BBA)\in\mathscr L^{\geq}$.
\hfill\QED
\end{corollary}
 
\section{Further applications of square roots}
\label{sec4}
First we express the generalized Schur complement of an operator
of $\mathscr L^{\geq}$ with the aid of its square root and derive 
a range description, cf.~\cite[Corollary~4 to Theorem~1]{Anderson}.
Let $\BBA\in\mathscr L^{\geq}$ and $(R,H)$ be a square root of $\BBA$.
Let $L$ be the orthogonal complement of $(RX)^{c}$ and $P$ be the 
orthoprojection onto $L$. Note that $L$ can be characterized by 
$L=\{h\in H\,{:}\, R^{\ast}h\in Y'\}$, which yields 
$R^{\ast}L=\ran R^{\ast}\cap Y'$.
\begin{proposition}
\label{prop4.1}
If $\BBA\in\mathscr L^{\geq}$, then $\mathscr S(\BBA)=R^{\ast}PR$.
\end{proposition}
\begin{pf}
Let $\begin{mmatrix} x \\ y \end{mmatrix}\in X\times Y$. An
application of \eqref{sec3:3} gives 
\[
\left\langle\mathscr S(\BBA)
\begin{pmatrix} x \\ y \end{pmatrix},
\begin{pmatrix} x \\ y \end{pmatrix}
\right\rangle =
\inf_{z\in X} \left\|
R\begin{pmatrix} x \\ y \end{pmatrix} -
R\begin{pmatrix} z \\ 0 \end{pmatrix}
\right\|^{2},
\]
which shows that 
$\left\langle\mathscr S(\BBA)
\begin{mmatrix} x \\ y \end{mmatrix},
\begin{mmatrix} x \\ y \end{mmatrix}
\right\rangle$ is the squared distance of 
$R\begin{mmatrix} x \\ y \end{mmatrix}$
to $RX$. Therefore, 
\[
\left\langle\mathscr S(\BBA)
\begin{pmatrix} x \\ y \end{pmatrix},
\begin{pmatrix} x \\ y \end{pmatrix}
\right\rangle =
\left\|PR
\begin{pmatrix} x \\ y \end{pmatrix} \right\|^{2} = 
\left\langle
R^{\ast}PR 
\begin{pmatrix} x \\ y \end{pmatrix},
\begin{pmatrix} x \\ y \end{pmatrix}
\right\rangle
\]
and the assertion follows from the polarization identity.
\hfill\QED
\end{pf}
\begin{proposition}
\label{prop4.2}
If $\BBA\in\mathscr L^{\geq}$ and $(R,H)$ and $(S,G)$ are
square roots of $\BBA$ and $\mathscr S(\BBA)$, resp., then
$\ran S^{\ast} = \ran R^{\ast}\cap Y'$.
\end{proposition}
\begin{pf}
Setting $R_{X}:=R\harp_{X}$ and $R_{Y}:= R\harp_{Y}$, we get 
\[
\BBA= 
\begin{pmatrix} R^{\ast}_{X} \\ R^{\ast}_{Y} \end{pmatrix}
(R_{X}\ R_{Y}) =
\begin{pmatrix} R^{\ast}_{X}R_{X}& R^{\ast}_{X} R_{Y} 
\\ R^{\ast}_{Y} R_{X} & R^{\ast}_{Y} R_{Y} 
\end{pmatrix},
\]
hence, 
\[
\sigma(\BBA)=R^{\ast}_{Y}R_{Y}-(R^{\ast[-1]}_{X}R^{\ast}_{X}R_{Y})^{\ast}
R^{\ast[-1]}_{X}R^{\ast}_{X}R_{Y}= R^{\ast}_{Y}P R_{Y}
\]
since $R^{\ast[-1]}_{X}R^{\ast}_{X}=I-P$. Thus 
$\mathscr S(\BBA)=R^{\ast}PR$ and $(PR,H)$ is a square root of 
$\mathscr S(\BBA)$. 
If $\begin{mmatrix} 0 \\ y'\end{mmatrix} \in X'\times Y'$ is
such that $R^{\ast}h = \begin{mmatrix} 0 \\ y'\end{mmatrix}$
for some $h\in H$, then $R^{\ast}_{X}h=0$, hence, $Ph=h$ and
$(PR)^{\ast}h=R^{\ast}h
=\begin{mmatrix} 0 \\ y'\end{mmatrix}$,
which implies that $\ran R^{\ast}\cap Y'\subseteq\ran(PR)^{\ast}
=\ran S^{\ast}$ by Corollary~\ref{cor2.6} (iii).
Since, obviously, $\ran S^{\ast}\subseteq Y'$ and 
$\ran S^{\ast}\subseteq\ran R^{\ast}$ by Corollaries~\ref{cor3.8} (i)
and \ref{cor2.6} (i), the assertion is proved.
\hfill\QED
\end{pf}
Our next result is a generalization of the Crabtree-Haynsworth 
quotient formula \cite{Crabtree}. To give it a nice form let us denote 
$\sigma(\BBA)=:{\BBA}/{A}$.
\begin{proposition}
\label{prop4.3}
Let $X$, $Y$, and $Z$ be linear spaces, 
\[
\BBD := 
\begin{pmatrix} A & B & B_{X} \\ 
B^{\sim} & D & B_{Y} \\
B^{\sim}_{X} & B^{\sim}_{Y} & D_{1}
\end{pmatrix}
\in \mathscr L^{\geq}(X\times Y\times Z, X'\times Y'\times Z'),
\]
and $\BBA := \begin{mmatrix} A & B \\ B^{\sim} & D \end{mmatrix}$. 
The operator ${\BBA}/{A}$ is the left upper corner of
${\BBD}/{A}$ and
\[
\quotient{{\BBD}/{A}}{{\BBA}/{A}}= {\BBD}/{\BBA}.
\]
\end{proposition}
\begin{pf}
Let $(R,H)$ be a minimal square root of $\BBA$,
$R_{X}:=R\harp_{X}$, $R_{Y}:=R\harp_{Y}$,
\[
E:= (R^{\ast})^{-1}
\begin{pmatrix} B_{X} \\ B_{Y}
\end{pmatrix},
\]
hence, $R^{\ast}_{X}E=B_{X}$,  $R^{\ast}_{Y}E=B_{Y}$. From
$R^{\ast[-1]}_{X}R^{\ast}_{X}=I-P$ one obtains
\begin{align*}
\BBD/A & = 
\begin{pmatrix}
R^{\ast}_{Y}R_{Y} & R^{\ast}_{Y} E \\
E^{\ast}R_{Y} & D_{1}
\end{pmatrix}
-\bigl(R^{\ast[-1]}_{X}
(R^{\ast}_{X}R_{Y},\ R^{\ast}_{X}E)
\bigr)^{\ast}
R^{\ast[-1]}_{X}
(R^{\ast}_{X}R_{Y},\ R^{\ast}_{X}E )\\
& =
\begin{pmatrix}
R^{\ast}_{Y}PR_{Y} & R^{\ast}_{Y}P E \\
E^{\ast}P R_{Y} & D_{1}-E^{\ast}(I-P)E
\end{pmatrix}
\end{align*}
and
\[
\BBA/A = R^{\ast}_{Y}R_{Y} - 
\bigl(R^{\ast[-1]}_{X}R^{\ast}_{X}R_{Y}\bigr)^{\ast}
R^{\ast[-1]}_{X}R^{\ast}_{X}R_{Y}
= R^{\ast}_{Y}PR_{Y},
\]
which shows that $\BBA/A$ is the left upper corner of
$\BBD/A$. Since $(PR_{Y},H)$ is a square root of $\BBA/A$,
one can compute 
\begin{align*}
\quotient{{\BBD}/{A}}{{\BBA}/{A}} & = 
D_{1} - E^{\ast}(I-P)E - 
\bigl((PR_{Y})^{\ast[-1]}R^{\ast}_{Y}PE\bigr)^{\ast}
(PR_{Y})^{\ast[-1]}R^{\ast}_{Y}PE \\
& = D_{1} -  E^{\ast}(I-P)E -E^{\ast}QE,
\end{align*}
where $Q$ denotes the orthoprojection onto $(\ran PR_{Y})^{c}$.
Comparing this with 
\[
\BBD/A = D_{1} - 
\bigl((R^{\ast})^{-1}R^{\ast}E\bigr)^{\ast}(R^{\ast})^{-1}R^{\ast}E
=D_{1}-E^{\ast}E,
\]
we can conclude that the assertion will be proved if we can 
show that the restriction of $I-P+Q$ to $\ran R$ is the 
identity. If $h\in\ran R$, then 
\[
h=R_{X}x+R_{Y}y = R_{X}x +(I-P)R_{Y}y+PR_{Y}y
\]
for some $\begin{mmatrix} x\\y\end{mmatrix}\in X\times Y$. Since
$R_{X}x+(I-P)R_{Y}y\in(\ran R_{X})^{c}$, there exists a sequence
$\{x_{n}\}_{n\in\BBN}$ of elements of $X$ such that
$\lim_{n\to\infty}R_{X}x_{n}=R_{X}x +(I-P)R_{Y}y$.
For $h_{n}:=R_{X}x_{n}+PR_{Y}y$, we have 
\[
(I-P+Q)h_{n}=(I-P+Q)(R_{X}x_{n}+PR_{Y}y)=R_{X}x_{n}+PR_{Y}y=h_{n}
\]
and therefore
\begin{align}
\notag
(I-P+Q)h & = \lim_{n\to\infty} (I-P+Q)h_{n} \\
\tag*{\QED}
& =
\lim_{n\to\infty} (R_{X}x_{n}+PR_{Y}y) = R_{X}x+R_{Y}y = h.
\end{align}
\end{pf}
We conclude this section with a criterion for non-negativity 
of operators of $\mathscr L^{h}$.
\begin{proposition}
\label{prop4.4}
Let 
$ \BBA = \begin{mmatrix} A & B \\ B^{\sim} & D\end{mmatrix}
\in\mathscr L^{h}$.
The operator $\BBA$ is non-negative if and only if the 
following two conditions are satisfied:
\begin{enumerate}
\item[(i)]
The operators $A$ and $D$ are non-negative.
\item[(ii)]
For any square roots $(R_{A},H_{A})$ and $(R_{D},H_{D})$ of
$A$ and $D$, resp., there exists a contraction 
$K\in\mathscr L(H_{D},H_{A})$ such that $B=R^{\ast}_{A}KR$
and $\ran K\subseteq (\ran R_{A})^{c}$.
\end{enumerate} 
\end{proposition}
\begin{pf}
If $\BBA$ is non-negative, assertion (i) is trivial. To prove 
(ii) let $(R,H)$ be a square root of $\BBA$ and $R_{X}:=R\harp_{X}$,
$R_{Y}=R\harp_{Y}$, hence 
\[
\BBA=
\begin{pmatrix} R^{\ast}_{X}R_{X} &  R^{\ast}_{X}R_{Y} 
\\ R^{\ast}_{Y}R_{X} & R^{\ast}_{Y}R_{Y} 
\end{pmatrix}.
\]
Let $(S_{A},G_{A})$ and $(S_{D},G_{D})$ be minimal square roots
of $A$ and $D$, resp.. According to Lemma~\ref{lemma2.8} there 
exist isometries $U_{A}\in\mathscr L(G_{A},H_{A})$, 
$V_{A}\in\mathscr L(G_{A},H)$, $U_{D}\in\mathscr L(G_{D},H_{D})$,
$V_{D}\in\mathscr L(G_{D},H)$
satisfying $U_{A}S_{A}=R_{A}$, $V_{A}S_{A}=R_{X}$, $U_{D}S_{D}=R_{D}$,
$V_{D}S_{D}=R_{Y}$. It follows 
\[
B = R^{\ast}_{X} R_{Y} =
R^{\ast}_{A}U_{A}V^{\ast}_{A}V_{D}U^{\ast}_{D}R_{D}
=R^{\ast}_{A} K R_{D},
\]
where $K:=U_{A}V^{\ast}_{A}V_{D}U^{\ast}_{D}\in\mathscr L
(H_{D},H_{A})$ is a contraction with 
$\ran K\subseteq (\ran R_{A})^{c}$.

Conversely, if (i) and (ii) are satisfied, then
\[
\BBA=
\begin{pmatrix} 
R^{\ast}_{A}R_{A} &  R^{\ast}_{A}KR_{D} \\ 
(R^{\ast}_{A}KR_{D})^{\sim} & R^{\ast}_{D}R_{D} 
\end{pmatrix}.
\]
Since 
\[
(R^{\ast}_{A}K R_{D})^{\sim} = (R'_{A}KR_{D})'\harp_{X} =
R'_{D}K' R^{\prime\prime}_{A}\harp_{X} =
R^{\ast}_{D} K^{\ast}R_{A},
\]
one obtains 
\[
\BBA=
\begin{pmatrix} 
R^{\ast}_{A} &  0 \\ 0 & R^{\ast}_{D} 
\end{pmatrix}
\begin{pmatrix} 
I  &  K \\ K^{\ast} & I  
\end{pmatrix}
\begin{pmatrix} 
R_{A} &  0 \\ 0 & R_{D} 
\end{pmatrix},
\]
which implies that $\BBA$ is non-negative.
\hfill\QED
\end{pf}
\section{Albert's theorem}
\label{sec5}

An application of Proposition~\ref{prop4.4} leads to a generalization
of an important criterion for non-negativity \cite{Albert}, which is often called
Albert's theorem in matrix theory. It should be
mentioned that Shmulyan~\cite[Theorem~1.7]{Shmulyan} had proved a
similar assertion even for bounded operators in Hilbert spaces ten
years earlier. We also mention the papers \cite{Ando} and
\cite{Hoeschel}.

\begin{theorem}
\label{theorem5.1}
An operator $\BBA=\begin{mmatrix} A & B \\ B^{\sim} & D\end{mmatrix}
\in\mathscr L^{h}$ is non-negative if and only if it is of positive
type and $\sigma(\BBA)$ is non-negative.
\end{theorem}  

\begin{pf}
If $\BBA\in\mathscr L^{\geq}$, Proposition~\ref{prop4.4} implies
that $\BBA$ is of positive type and $R^{\ast[-1]}_{A}B=KR_{D}$
for some contraction $K\in\mathscr L(H_{D},H_{A})$. It follows
\[
\sigma(\BBA) = R^{\ast}_{D}R_{D}-(KR_{D})^{\ast}KR_{D} =
R^{\ast}_{D}(I-K^{\ast}K)R_{D}\geq 0
\]
Conversely, let $\BBA\in\mathscr L^{+}$ and $\sigma(\BBA)\in\mathscr
L^{\geq}(Y,Y')$. If $(R_{A},H_{A})$ and $(R_{D},H_{D})$ are square
roots of $A$ and $D$, resp., one has
\[
\|R^{\ast[-1]}_{A}By\|^{2} =
\langle\omega(A,B)y,y\rangle \leq\langle Dy,y\rangle =
\|R_{D}y\|^{2},\quad y\in Y,
\]  
which yields $KR_{D}=R^{\ast[-1]}B$, hence, $R^{\ast}_{A}KR_{D}=B$
for some contraction $K\in\mathscr L(H_{D},H_{A})$. An application
of Proposition~\ref{prop4.4} completes the proof.
\hfill\QED
\end{pf}

The preceding theorem can be used to study the set $\mathscr
L^{\geq}$ as well as the set $\mathscr L^{+}$ and to establish
interrelations between these two sets. A first result is the
inclusion $\mathscr L^{\geq}\subseteq\mathscr L^{+}$.
For a positive pair $(A,B)$ set
\[
\BBA_{ex}:=
\begin{pmatrix}A & B \\ B^{\sim} & \omega(A,B) \end{pmatrix}
\in\mathscr L^{+}.
\]

\begin{corollary}
\label{cor5.2}
Two operators $A\in\mathscr L(X,X')$ and $B\in\mathscr L(Y,X')$
form a positive pair if and only if the set
\[
  \mathscr A := \left \{ \BBA\in\mathscr L^{\geq}\;{:}\;
  \BBA=\begin{mmatrix} A & B \\ B^{\sim} & D\end{mmatrix}
  \text{ for some $D\in\mathscr L^{\geq}(Y,Y')$}
  \right\}
\]
is non-empty. If $(A,B)$ is a positive pair, the operator $\BBA_{ex}$
is the minimal element of $\mathcal A$.
\hfill\QED
\end{corollary}

\begin{corollary}
\label{cor5.3}
If $\BBA\in\mathscr L^{\geq}$, the set 
\[
\mathscr A_{1} :=
\bigl\{\BBA_{1}\in\mathscr L^{\geq}\;{:}\;\BBA_{1}\leq \BBA
\text{ and } X\subseteq \ker \BBA_{1}\bigr\}
\]
is non-empty and $\mathscr S(\BBA)$ is its maximal element.
\end{corollary}

\begin{pf}
Corollaries~3.8 (i) and 3.9 imply that $\mathscr S(\BBA)\in \mathscr
A_{1}$. If $\BBA = \begin{mmatrix} A & B \\ B^{\sim} & D\end{mmatrix}$
and $\BBA_{1}\in\mathscr A_{1}$, then $\BBA_{1}$ has representation
$\BBA_{1} = \begin{mmatrix} 0 & 0 \\ 0 & D_{1}\end{mmatrix}$ and
$D-\omega(A,B)-D_{1}\geq 0$, hence, $\BBA_{1}\leq\mathscr S(\BBA)$
by Theorem~\ref{theorem5.1}.
\hfill\QED
\end{pf}

\begin{corollary}
\label{cor5.4}
Let $\BBA =  \begin{mmatrix} A & B \\ B^{\sim} & D\end{mmatrix}
\in\mathscr L^{h}$. The operator $\BBA$ belongs to $\mathscr L^{+}$
if and only if there exists an operator $\BBA_{1}\in \mathscr
L^{h}$ satisfying $X\subseteq\ker \BBA_{1}$ and $\BBA_{1}\leq\BBA$.
\end{corollary}

\begin{pf}
  If $\BBA\in\mathscr L^{+}$, the operator $\BBA_{1}:=\mathscr S(\BBA)$
  has all properties claimed. Conversely, if there exists an operator
  $\BBA_{1}$ satisfying all conditions, it has the form
  $\BBA_{1}=\begin{mmatrix}0 & 0 \\ 0 & D_{1}\end{mmatrix}$, where
  $D_{1}\in\mathscr L(Y,Y')$ and $\BBA-\BBA_{1}=
  \begin{mmatrix} A & B \\ B^{\sim} & D-D_{1}\end{mmatrix}\in
  \mathscr L^{\geq}$. It follows from Theorem~\ref{theorem5.1} that
  $(A,B)$ is a positive pair, hence, $\BBA\in\mathscr L^{+}$.
\hfill\QED
\end{pf}

Another application of Theorem~\ref{theorem5.1} gives an expression of
the supremum occuring in Lemma~\ref{lemma3.2}.

\begin{corollary}
\label{cor5.5}
If $(A,B)$ is a positive pair, then
\begin{equation}
\label{sec5:1}
\sup_{x\in X} \frac{|\langle By,x\rangle|^{2}}{\langle Ax,x\rangle}
=\langle\omega(A,B)y,y\rangle,\quad y\in Y.
\end{equation}  
\end{corollary}

\begin{pf}
  Let $y\in Y$. Since $\BBA_{ex}\in\mathscr L^{\geq}$ by
  Corollary~\ref{cor5.2}, one has
  \[
    |\langle By,x\rangle|^{2}\leq\langle Ax,x\rangle
    \langle\omega(A,B)y,y\rangle,
  \]
  which yields
  \[
    \frac{|\langle By,x\rangle|^{2}}{\langle Ax,x\rangle}
    \leq \langle\omega(A,B)y,y\rangle, \quad x\in X,
  \]
  if one takes into account the convention $\frac{0}{0}:= 0$. Thus,
  \eqref{sec5:1} has been proved if $Ty = R^{\ast[-1]}By=0$, where $(R,H)$ is
  a minimal square root of $A$. Now assume that $Ty\neq0$.
  There exists a sequence $\{x_{n}\}_{n\in\BBN}$ of elements of $X$
  such that $Rx_{n}\neq0$, $n\in \BBN$, and $\lim_{n\to\infty}Rx_{n}=
  Ty$. It follows
  \begin{align}
    \notag
    \lim_{n\to\infty} 
    \frac{|\langle By,x_{n}\rangle|^{2}}{\langle Ax_{n},x_{n}\rangle}
    & =
    \lim_{n\to\infty}
    \frac{|\langle R^{\ast}Ty,x_{n}\rangle|^{2}}%
    {\langle R^{\ast}Rx_{n},x_{n}\rangle}
    =
    \lim_{n\to\infty}
      \frac{|(Ty\vertic Rx_{n})|^{2}}{\|Rx_{n}\|^{2}} \\
    \tag*{\QED}
    & =
    \|Ty\|^{2} =  \langle\omega(A,B)y,y\rangle.
  \end{align}  
\end{pf}

\begin{corollary}
\label{cor5.6}
Let $(A_{j},B_{j})$ with $A_{j}\in\mathscr(X,X')$,
$B_{j}\in\mathscr(Y,X')$, $j=1,2$, be positive pairs. Then
$(A_{1}+A_{2},B_{1}+B_{2})$ is a positive pair and
\begin{equation}
\label{sec5:2}
  \omega(A_{1}+A_{2},B_{1}+B_{2})\leq
  \omega(A_{1},B_{1})+\omega(A_{2},B_{2}).
\end{equation}
\end{corollary}

\begin{pf}
Since the operators $(\BBA_{j})_{ex}$, $j=1,2$ are non-negative, it
follows
\[
\begin{pmatrix}
  A_{1}+A_{2} & B_{1}+ B_{2} \\
  B^{\sim}_{1}+B^{\sim}_{2} & \omega(A_{1},B_{1})+\omega(A_{2},B_{2})
\end{pmatrix}
\in\mathscr L^{\geq},
\]
hence, \eqref{sec5:2} by Corollary~\ref{cor5.2}.
\hfill\QED
\end{pf}

\begin{corollary}
\label{cor5.7}
If $\BBA_{j}\in\mathscr L^{+}$, $j=1,2$, then
$\BBA_{1}+\BBA_{2}\in\mathscr L^{+}$ and
\begin{equation}
  \tag*{\QED}
  \mathscr S(\BBA_{1})+\mathscr S(\BBA_{2})
  \leq \mathscr S(\BBA_{1}+\BBA_{2}).
\end{equation}
\end{corollary}

A subset $\mathscr A$ of $\mathscr L^{h}$ is called bounded below if
there exists $\BBA_{1}\in\mathscr L^{h}$ such that $\BBA_{1}\leq\BBA$
for all $\BBA\in\mathscr A$. An operator $\BBA_{0}\in\mathscr L^{h}$
is called an infimum of $\mathscr A$ if the following conditions are
satisfied:
\begin{enumerate}
\item[a)] $\BBA_{0}\leq\BBA$ for all $\BBA\in\mathscr A$,
\item[b)] $\BBA_{1}\leq\BBA_{0}$ for all $\BBA_{1}\in\mathscr L^{h}$
  such that $\BBA_{1}\leq\BBA$, $\BBA\in\mathscr A$.
\end{enumerate}
If an infimum of $\mathscr A$ exists, it is unique. Recall that any
set $\mathscr A$, which is bounded from below and directed downwards
(i.e. for all $\BBA_{1},\BBA_{2}\in\mathscr A$ there exists
$\BBA\in\mathscr A$ such that $\BBA\leq\BBA_{1}$ and
$\BBA\leq\BBA_{2}$), possesses an infimum. Particularly, if
$\{\BBA_{n}\}_{n\in\BBN}$ is a decreasing sequence of operators of
$\mathscr L^{h}$, which is bounded from below, there exists an infimum
$\BBA_{0}$ and $\langle\BBA_{0}z_{1},z_{2}\rangle=\lim_{n\to\infty}
\langle\BBA_{n} z_{1},z_{2}\rangle$ for all $z_{1},z_{2}\in X\times Y$.

\begin{corollary}
\label{cor5.8}
Let $\mathscr A$ be a subset of $\mathscr L^{h}$, which has an infimum
$\BBA_{0}$. The operator $\BBA_{0}$ belongs to $\mathscr L^{+}$ if and
only if the set $\mathcal S(\mathscr A):=\{\mathscr S(\BBA)\;{:}\;\BBA\in\mathscr A\}$
is bounded from below. In this case $\mathscr S(\BBA_{0})$ is the
infimum of $\mathcal S(\mathscr A)$.
\end{corollary}

\begin{pf}
If $\BBA_{0}\in\mathscr L^{+}$, the set $\mathcal S(\mathscr A)$ is
bounded from below since $\mathscr S(\BBA_{0})\leq\mathscr S(\BBA)$,
$\BBA\in\mathscr A$, by Corollary~\ref{cor3.8} (ii). Conversely,
assume
that there exists $\BBA_{1}\in\mathscr L^{h}$ such that
$\BBA_{1}\leq\mathscr S(\BBA)$ for all $\BBA\in\mathscr A$. It follows
$-\BBA_{1}\geq-\mathscr S(\BBA)$, which yields $-\BBA_{1}\in\mathscr
L^{+}$
by Corollary~\ref{cor5.4} and $\mathscr S(-\BBA_{1})\geq\mathscr
S(-\mathscr S(\BBA))=-\mathscr S(\BBA)$, hence, $-\mathscr
S(-\BBA_{1})\leq \mathscr S(\BBA)\leq\BBA$, $\BBA\in\mathscr A$, 
by Corollary~\ref{cor3.8}. One obtains $-\mathscr S(-\BBA_{1})\leq
\BBA_{0}$ and therefore $\BBA_{0}\in\mathscr L^{+}$ by
Corollary~\ref{cor5.4}. Moreover, $\mathscr S(\BBA_{0})\leq
\mathscr S(\BBA)$, $\BBA\in\mathscr A$, and
\[
  \BBA_{1}=-(-\BBA_{1})\leq-\mathscr S(-\BBA_{1})=
  \mathscr S\bigl(-\mathscr S(-\BBA_{1})\bigr)\leq\mathscr S(\BBA_{0})
\]  
by Corollary~\ref{cor3.8}, which implies that $\mathscr S(\BBA_{0})$
is the infimum of $\mathcal S(\mathscr A)$.
\hfill\QED
\end{pf}

\section{Extremal operators}
\label{sec6}

An operator $\BBA\in\mathscr L^{+}$ was called an extremal operator by
M.G. Kre\u{\i}n \cite{Krein} if $\mathscr S(\BBA)=0$. Since $\BBA=\mathscr
S(\BBA)+\BBA_{ex}$, an operator is extremal if and only if it has the
form
\[
  \BBA=\BBA_{ex} =
  \begin{pmatrix}A & B \\ B^{\sim} & \omega(A,B)\end{pmatrix}
\]  
for some positive pair $(A,B)$. Particularly, any extremal operator is
non-negative. Applying Proposition~\ref{prop3.7} we can give several
criteria for an operator to be extremal.

\begin{lemma}
\label{lemma6.1}
Let $\BBA\in\mathscr L^{\geq}$. The following assertions are
equivalent:
\begin{enumerate}
\item[(i)] The operator is extremal.
\item[(ii)] For all $\begin{mmatrix}x\\y\end{mmatrix}\in X\times Y$
  and arbitrary $\varepsilon>0$ there exists $z\in X$ such that
  \[
    \left\langle\BBA\begin{mmatrix}x-z\\y\end{mmatrix},
      \begin{mmatrix}x-z\\y\end{mmatrix}\right\rangle <\varepsilon.
  \]  
\item[(iii)] For any square root $(R,H)$ of $\BBA$ the spaces
  $(RX)^{c}$ and $(\ran R)^{c}$ coincide.
\item[(iv)] For any square root $(R,H)$ of $\BBA$ one has
  $\ran R^{\ast}\cap Y'=\{0\}$.
\end{enumerate}
\end{lemma}  

\begin{pf}
 The equivalence of (i) and (ii) is an immediate
 consequence of \eqref{sec3:3}. To prove (i)  $\Leftrightarrow$ (iii),
 choose a minimal square root $(R,H)$ of $\BBA$ and let $L$ and $P$
 be defined as in Proposition~\ref{prop4.1}. Then $\mathscr
 S(\BBA)=R^{\ast}PR=0$ if and only if $P=0$ or, equivalently,
 $L=\{0\}$, which in turn is equivalent to $(RX)^{c}=(\ran R)^{c}$.
 The equivalence of (iii) and (iv) follows from the equality
 $R^{\ast}L=\ran R^{\ast}\cap Y'$.
 \hfill\QED 
\end{pf}

Let $(A,B)$ be a positive pair and $(R,H)$ a square root of $A$. 
Recall the notation (4.1) of the operator $T:=R^{\ast[-1]}B$. Moreover, 
let $P_{B}$ be the orthoprojection onto
$(\ran T)^{c}$. Since the operator $\BBA_{ex}$ is
non-negative, from Theorem~\ref{theorem5.1} one can conclude that
$(\omega(A,B),B^{\sim})$ is a positive pair as well changing the roles
of $X$ and $Y$. Thus, the
operators $T^{\ast[-1]}B^{\sim}$ and
\[
  \omega\bigl(\omega(A,B),B^{\sim}\bigr)=
  \bigl[T^{\ast[-1]}B^{\sim}\bigr]^{\ast}
  T^{\ast[-1]}B^{\sim}
\]
can be defined.  

\begin{lemma}
\label{lemma6.2}
The equalities $T^{\ast[-1]}B^{\sim} = P_{B}R$
and  $\omega\bigl(\omega(A,B),B^{\sim}\bigr) =R^{\ast}P_{B}R$
hold true.
\end{lemma}

\begin{pf}
 The second equality is an immediate consequence of the first one.
 To prove the first equality we shall show that
\begin{equation}
\label{sec6:1}
\bigl( T^{\ast[-1]}B^{\sim}x\vertic h\bigr)
= (P_{B}R x\vertic h)\quad
\text{for all $x\in X$ and $h\in H$.}
\end{equation}
Since $\ran T^{\ast[-1]}B^{\sim}
\subseteq (\ran T)^{c}$ it is enough to prove
\eqref{sec6:1} for $x\in X$ and $h\in\ran T$.
If $h= Ty$ for some $y\in Y$, we get
\begin{align*}
\bigl( T^{\ast[-1]}B^{\sim}x\vertic h\bigr)
&  =
\bigl( T^{\ast[-1]}B^{\sim}x\vertic
  Ty\bigr) \\
& =
  \bigl\langle T^{\ast}
  T^{\ast[-1]} B^{\sim}x,y \bigr\rangle
  = \langle B^{\sim} x,y\rangle
\end{align*}
and
\begin{align*}
  (P_{B}Rx\vertic h) & = (P_{B}Rx\vertic Ty) = (Rx\vertic
  Ty) \\
  & = 
  \overline{\langle R^{\ast}Ty,x\rangle} =
  \overline{\langle By,x\rangle} =\langle B^{\sim}x,y\rangle,
\end{align*}
hence, \eqref{sec6:1}.
\hfill\QED
\end{pf}

{}From Corollary~\ref{cor5.2} it follows that
$\omega\bigl(\omega(A,B),B^{\sim}\bigr)$ is a minimal element of the
set
\[
  \left\{A_{1}\in\mathscr L(X,X')\;{:}\;
    \begin{mmatrix} A_{1} & B \\ B^{\sim} & \omega(A,B)
    \end{mmatrix}\in\mathscr L^{\geq}\right\}  
\]
Note also that
\[
\omega\Bigl(\omega\bigl(\omega(A,B),B^{\sim}\bigr),B\Bigr)
= \omega(A,B),
\]
cf. \cite[Proposition~1.4 (A)]{Niemiec}. We call an extremal operator
$\BBA_{ex}=\begin{mmatrix}A&B\\B^{\sim}& \omega(A,B)\end{mmatrix}$
doubly extremal if $\omega\bigl(\omega(A,B),B^{\sim}\bigr)=A$.

In the case of bounded operators on Hilbert spaces the remaining
results of the present section were proved by Pekarev and Shmulyan \cite{Pekarev}
and partly rediscovered by Niemiec \cite{Niemiec}. We mention that
Niemiec's proofs are based on Douglas' theorem and do not make explicit
use of $2\times 2$~block operators.

\begin{proposition}
\label{prop6.3}
An operator $\BBA_{ex}$ is doubly extremal if and only if 
\begin{equation}
\label{sec6:2}
\ker T^{\ast} = \ker R^{\ast}
\end{equation}
for any square root $(R,H)$ of $A$.
\end{proposition}

\begin{pf}
 According to Lemma~\ref{lemma6.2}, $\BBA$ is doubly extremal if and
 only if $R^{\ast}P_{B}R=R^{\ast}R$. If
 $\ker T^{\ast} = \ker R^{\ast}$ or, equivalently,
 $(\ran T)^{c}=(\ran R)^{c}$, it follows $P_{B}R=R$,
 hence, $R^{\ast}P_{B}R=R^{\ast}R$. Conversely, assume
 $R^{\ast}P_{B}R=R^{\ast}R$, which yields $\|P_{B}Rx\|=\|Rx\|$,
 $x\in X$, hence, $(\ran R)^{c}\subseteq\ran P_{B}$ and $\ker
 P_{B}\subseteq\ker R^{\ast}$. Since $\ran P_{B}=(\ran
 T)^{c}$ or $\ker P_{B}= \ker T^{\ast}$,
we get 
\begin{equation}
\label{sec6:3}
\ker T^{\ast} \subseteq \ker R^{\ast}.
\end{equation}
On the other hand, if $h\in\ker R^{\ast}$, then $h$ is orthogonal to
$(\ran R)^{c}$ and
\[
  0 = \bigl(h\vertic Ty \bigr)
    = \bigl\langle T^{\ast}h,y\bigr\rangle,
  \quad y\in Y,
\]  
thus, $\ker R^{\ast}\subseteq \ker T^{\ast}$.
Taking into account \eqref{sec6:3} we obtain the desired equality.
\hfill\QED
\end{pf}

The preceding assertion shows that the equality \eqref{sec6:2}
does not depend on the choice of the square root $(R,H)$ of $A$ and
that in the case of a minimal square root the operator $\BBA_{ex}$
is doubly extremal if and only if $\ker T^{\ast}=\{0\}$.
Moreover, writing \eqref{sec6:2} in the equivalent form
$(\ran T)^{c}=(\ran R)^{c}$ we obtain a generalization of
\cite[Theorem~1.6]{Pekarev}. It means that $\BBA_{ex}$ is doubly
extremal if and only if the inverse image of $\ran B$ under the map
$R^{\ast}$ is dense in $H$.

\begin{corollary}
\label{cor6.4}
If $\ran R^{\ast}=\ran B$, the operator $\BBA_{ex}$ is doubly extremal.
\hfill\QED 
\end{corollary}

To give another criterion for $\BBA_{ex}$ to be doubly extremal we
equip the space $Y'$ with $\sigma(Y',Y)$-topology, i.e. the smallest
topology such that for arbitrary $y\in Y$, the functional $y'\mapsto
\langle y',y\rangle$ is continuous on $Y'$. Let $(A,B)$ be a positive
pair and $(R,H)$ a square root of $A$. Denote by $H_{1}$ the subspace of
all $h\in H$ such that there exists a sequence $\{x_{n}\}_{n\in\BBN}$
of elements of $X$ with the following properties:
\begin{enumerate}
\item[a)]
  $\lim_{n\to\infty}Rx_{n}=h$ with respect to the norm topology of
  $H$,
\item[b)]
  $\lim_{n\to\infty} B^{\sim}x_{n}=0$ with respect to the
  $\sigma(Y',Y)$-topology.
\end{enumerate}

\begin{lemma}
\label{lemma6.5}
The space $H_{1}$ is equal to $(\ran R)^{c}\cap\ker T^{\ast}$.
\end{lemma}  

\begin{pf}
An element $h\in H$ belongs to $H_{1}$ if and only if there exists a
sequence $\{x_{n}\}_{n\in\BBN}$ of elements of $X$ such that
$\lim_{n\to\infty}Rx_{n}=h$ and for all $y\in Y$,
\begin{align}
\notag  
\bigl\langle T^{\ast}h,y\bigr\rangle
& =
\lim_{n\to\infty} \bigl(Rx_{n}\vertic Ty\bigr) \\
\tag*{\QED}  
& =
\lim_{n\to\infty} (x_{n}\vertic By) =
\lim_{n\to\infty}(B^{\sim}x_{n}\vertic y) = 0.
\end{align}      
\end{pf}

\begin{proposition}
\label{prop6.6}
An operator $\BBA_{ex}$ is doubly extremal if and only if $H_{1}=\{0\}$.
\end{proposition}

\begin{pf}
Since  $H_{1}=\{0\}$ if and only if
$\ker R^{\ast}=\ker T^{\ast}$
by Lemma~\ref{lemma6.5}, the assertion follows from
Proposition~\ref{prop6.3}.
\hfill\QED
\end{pf}  

\begin{corollary}
\label{cor6.7}
If an operator $\BBA_{ex}$ is doubly extremal, then $\ker A=\ker
B^{\sim}$. If $\ker A=\ker B^{\sim}$ and $\ran R$ is closed, then
$\BBA_{ex}$ is doubly extremal.  
\end{corollary}

\begin{pf}
 If $\BBA_{ex}$ is doubly extremal, then
$\bigl (T^{\ast[-1]}B^{\sim},H\bigr)$
is a square root of $A$, hence,
$\ker B^{\sim}\subseteq \ker T^{\ast[-1]}B^{\sim}=\ker A$
by Lemma~\ref{lemma2.4}. The first assertion of the corollary follows since
$\ker A\subseteq \ker B^{\sim}$ by Lemma~\ref{lemma3.2}. Now
assume that  $\ker A=\ker B^{\sim}$ and $\ran R$ is closed. If
$h\in H_{1}$, there exist $x\in X$ and a sequence
$\{x_{n}\}_{n\in\BBN}$ of elements of $X$ such that
$h=Rx=\lim_{n\to\infty}Rx_{n}$ and
$\lim_{n\to\infty}\langle B^{\sim}x_{n},y\rangle = 0$
for all $y\in Y$. It follows
\begin{align*}
  \langle B^{\sim}x,y\rangle
  & =
  \overline{\bigl\langle R^{\ast}Ty,x\bigr\rangle}
   =
  \bigl(Rx\vertic Ty\bigr) \\
  & =
  \lim_{n\to\infty} \bigl(Rx_{n}\vertic Ty\bigr)
  =
  \lim_{n\to\infty} \overline{\langle By,x_{n}\rangle} \\
  & =
  \lim_{n\to\infty} \langle B^{\sim}x_{n},y\rangle = 0, \quad y\in Y,
\end{align*}
which implies that $x\in\ker B^{\sim}=\ker A=\ker R$ and
$h=0$. An application of Proposition~\ref{prop6.6} completes
the proof.
\hfill\QED
\end{pf}




\bibliographystyle{amsplain}
\bibliography{schur}

\providecommand{\bysame}{\leavevmode\hbox to3em{\hrulefill}\thinspace}
\providecommand{\MR}{\relax\ifhmode\unskip\space\fi MR }
\providecommand{\MRhref}[2]{%
  \href{http://www.ams.org/mathscinet-getitem?mr=#1}{#2}
}
\providecommand{\href}[2]{#2}
\begin{thebibliography}{10}

\bibitem{Albert}
Arthur Albert, \emph{Conditions for positive and nonnegative definiteness in
  terms of pseudoinverses}, SIAM J. Appl. Math. \textbf{17} (1969), 434--440.
  \MR{0245582}

\bibitem{Anderson}
W.~N. Anderson, Jr. and G.~E. Trapp, \emph{Shorted operators. {II}}, SIAM J.
  Appl. Math. \textbf{28} (1975), 60--71. \MR{0356949}

\bibitem{Ando}
T.~And\^o, \emph{Truncated moment problems for operators}, Acta Sci. Math.
  (Szeged) \textbf{31} (1970), 319--334. \MR{0290157}

\bibitem{Crabtree}
Douglas~E. Crabtree and Emilie~V. Haynsworth, \emph{An identity for the {S}chur
  complement of a matrix}, Proc. Amer. Math. Soc. \textbf{22} (1969), 364--366.
  \MR{0255573}

\bibitem{Douglas}
R.~G. Douglas, \emph{On majorization, factorization, and range inclusion of
  operators on {H}ilbert space}, Proc. Amer. Math. Soc. \textbf{17} (1966),
  413--415. \MR{0203464}

\bibitem{FKKlotz}
Bernd Fritzsche, Bernd Kirstein, and Lutz Klotz, \emph{Completion of
  non-negative block operators in {B}anach spaces}, Positivity \textbf{3}
  (1999), no.~4, 389--397. \MR{1721580}

\bibitem{Gorniak84}
J.~G\'orniak, \emph{Locally convex spaces with factorization property}, Colloq.
  Math. \textbf{48} (1984), no.~1, 69--79. \MR{750756}

\bibitem{Gorniak80}
J.~G\'orniak and A.~Weron, \emph{Aronszajn-{K}olmogorov type theorems for
  positive definite kernels in locally convex spaces}, Studia Math. \textbf{69}
  (1980/81), no.~3, 235--246. \MR{647140}

\bibitem{Gorniak83}
Janusz G\'orniak, Andrzej Makagon, and Aleksander Weron, \emph{An explicit form
  of dilation theorems for semispectral measures}, Prediction theory and
  harmonic analysis, North-Holland, Amsterdam-New York, 1983, pp.~85--111.
  \MR{708520}

\bibitem{Hoeschel}
Hans-Peter H\"oschel, \emph{{\"U}ber die {P}seudoinverse eines zerlegten
  positiven linearen {O}perators}, Math. Nachr. \textbf{74} (1976), 167--172.
  \MR{0425637}

\bibitem{Krein}
M.~G. Kre\u{\i}n, \emph{The theory of self-adjoint extensions of semi-bounded
  {H}ermitian transformations and its applications. {I}}, Rec. Math. [Mat.
  Sbornik] N.S. \textbf{20(62)} (1947), 431--495. \MR{0024574}

\bibitem{Niemiec}
Piotr Niemiec, \emph{Generalized absolute values and polar decompositions of a
  bounded operator}, Integral Equations Operator Theory \textbf{71} (2011),
  no.~2, 151--160. \MR{2838139}

\bibitem{Pekarev}
\`E.~L. Pekarev and Yu.~L. Shmulyan, \emph{Parallel addition and parallel
  subtraction of operators}, Izv. Akad. Nauk SSSR Ser. Mat. \textbf{40} (1976),
  no.~2, 366--387, 470. \MR{0410429}

\bibitem{Pusz}
W.~Pusz and S.~L. Woronowicz, \emph{Functional calculus for sesquilinear forms
  and the purification map}, Rep. Mathematical Phys. \textbf{8} (1975), no.~2,
  159--170. \MR{0420302}

\bibitem{Schwartz}
Laurent Schwartz, \emph{Sous-espaces hilbertiens d'espaces vectoriels
  topologiques et noyaux associ\'es (noyaux reproduisants)}, J. Analyse Math.
  \textbf{13} (1964), 115--256. \MR{0179587}

\bibitem{Sebestyen}
Zolt\'an Sebesty\'en, \emph{Operator extensions on {H}ilbert space}, Acta Sci.
  Math. (Szeged) \textbf{57} (1993), no.~1-4, 233--248. \MR{1243281}

\bibitem{Shmulyan}
Yu.~L. Shmulyan, \emph{An operator {H}ellinger integral}, Mat. Sb. (N.S.)
  \textbf{49 (91)} (1959), 381--430. \MR{0121662}

\bibitem{Shmulyan1}
\bysame, \emph{Two-sided division in the ring of operators}, Mat. Zametki
  \textbf{1} (1967), 605--610. \MR{0217640}

\bibitem{Tarcsay}
Zs. Tarcsay, \emph{On the parallel sum of positive operators, forms, and
  functionals}, Acta Math. Hungar. \textbf{147} (2015), no.~2, 408--426.
  \MR{3420586}

\bibitem{Titkos}
T.~Titkos, \emph{On means of nonnegative sesquilinear forms}, Acta Math.
  Hungar. \textbf{143} (2014), no.~2, 515--533. \MR{3233548}

\bibitem{Vakhaniya}
N.~N. Vakhaniya, \emph{Probabilistic distributions in linear spaces}, Sakharth.
  SSR Mecn. Akad. Gamothvl. Centr. \v{S}rom. \textbf{10} (1971), no.~3, 155,
  Russian, English translation: North-Holland Publishing Co., New
  York-Amsterdam 1981. \MR{0296988}

\bibitem{Vakhaniya1}
N.~N. Vakhaniya, V.~I. Tarieladze, and S.~A. Chobanyan, \emph{Probability
  {D}istributions in banach {S}paces}, Nauka, Moscow, 1985, Russian, English
  translation: D.~Reidel Publishing Co, Dordrecht 1987. \MR{787803}

\bibitem{Vainberg}
M.~M. Va\u{\i}nberg, \emph{Variational method and method of monotone operators
  in the theory of nonlinear equations}, Nauka, Moscow, 1972, Russian, English
  Translation: Halsted Press, New York-Toronto 1973. \MR{0467427}

\bibitem{Zhang}
Fuzhen Zhang (ed.), \emph{The {S}chur {C}omplement and its {A}pplications},
  Numerical Methods and Algorithms, vol.~4, Springer-Verlag, New York, 2005.
  \MR{2160825}

\end{thebibliography}







\end{document}